\documentstyle[12pt,subeqn]{article}
\newfont{\twelvemsb}{msbm10 scaled\magstep1}
\newfont{\eightmsb}{msbm8}
\newfont{\sixmsb}{msbm6}
\newfam\msbfam
\textfont\msbfam=\twelvemsb
\scriptfont\msbfam=\eightmsb
\scriptscriptfont\msbfam=\sixmsb
\catcode`\@=11
\def\Bbb{\ifmmode\let\next\Bbb@\else
  \def\next{\errmessage{Use \string\Bbb\space only in math mode}}\fi\next}
\def\Bbb@#1{{\Bbb@@{#1}}}
\def\Bbb@@#1{\fam\msbfam#1}
\newfont{\twelvegoth}{eufm10 scaled\magstep1}
\newfont{\tengoth}{eufm10}
\newfont{\eightgoth}{eufm8}
\newfont{\sixgoth}{eufm6}
\newfam\gothfam
\textfont\gothfam=\twelvegoth
\scriptfont\gothfam=\eightgoth
\scriptscriptfont\gothfam=\sixgoth
\def\frak{\ifmmode\let\next\frak@\else
  \def\next{\errmessage{Use \string\frak\space only in math mode}}\fi\next}
\def\frak@#1{{\fam\gothfam{{#1}}}}
\catcode`@=12

\newcommand{\be}{\begin{equation}}
\newcommand{\ee}{\end{equation}}
\newcommand{\bea}{\begin{eqnarray}}
\newcommand{\ena}{\end{eqnarray}}
\newcommand{\sect}[1]{\setcounter{equation}{0}\section{#1}}

\newcommand{\cle}{\Big[}
\newcommand{\cri}{\Big]}
\newcommand{\ale}{\Big\{}
\newcommand{\ari}{\Big\}}
\newcommand{\sle}{\Big[\hspace{-3pt}\Big[}
\newcommand{\sri}{\Big]\hspace{-3pt}\Big]}
\newcommand{\ZZ}{{\Bbb Z}}
\newcommand{\CC}{{\Bbb C}}
\newcommand{\ZP}{{\ZZ'}}
\newcommand{\del}{\delta}
\newcommand{\eps}{\varepsilon}
\newcommand{\fra}{{\frak a}}
\newcommand{\frb}{{\frak b}}
\newcommand{\half}{\frac{1}{2}}
\newcommand{\smbox}[1]{\ \mbox{#1}\ }
\newcommand{\medbox}[1]{\quad\mbox{#1}\quad}
\newcommand{\bigbox}[1]{\qquad\mbox{#1}\qquad}
\newcommand{\cA}{{\cal A}}
\newcommand{\cE}{{\cal E}}
\newcommand{\cU}{{\cal U}}
\newcommand{\ta}{{\tilde a}}
\newcommand{\tA}{{\tilde A}}
\newcommand{\tK}{{\tilde K}}
\newcommand{\sa}{{\widehat A(M-1,N-1)}}
\newcommand{\uqa}{{\cU_q(\widehat A(M-1,N-1))}}
\newcommand{\uqsa}{{\cU_q(A(M-1,N-1))}}

\newtheorem{theorem}{Theorem}

\begin{document}
\newpage
\pagestyle{empty}
\setcounter{page}{0}
\vfill
\begin{center}
  
  {\Large {\bf ANYONIC REALIZATION OF THE 

  \vspace{5mm}

  QUANTUM AFFINE LIE SUPERALGEBRAS 
   
  \vspace{5mm}
  
  $\uqa$}}\\[1cm] 
  
  \vspace{7mm}
  
  {\large L. Frappat$^{a}$
          \footnote{On leave of absence from Laboratoire de Physique 
                    Th\'eorique ENSLAPP.},
    A. Sciarrino$^{b}$,
    S. Sciuto$^{c}$, 
    P. Sorba$^{d}$
  }

  \vspace{7mm}
  
  {\em $^{a}$ Centre de Recherches Math\'ematiques}
    \\
  {\em Universit\'e de Montr\'eal, Canada}

  \vspace{5mm}

  {\em $^{b}$ Dipartimento di Scienze Fisiche, Universit\`a di Napoli
       ``Federico II''}
     \\
  {\em and I.N.F.N., Sezione di Napoli, Italy}

  \vspace{5mm}

  {\em $^{c}$ Dipartimento di Fisica Teorica, Universit\`a di Torino}
     \\
  {\em and I.N.F.N., Sezione di Torino, Italy}

  \vspace{5mm}

  {\em $^{d}$ Laboratoire de Physique Th\'eorique ENSLAPP}
    \\
  {\em Annecy-le-Vieux et Lyon, France}

\end{center}

\vfill

\begin{abstract}
We give a realization of the quantum affine Lie superalgebras $\uqa$ in 
terms of anyons defined on a one or two-dimensional lattice, the deformation 
parameter $q$ being related to the statistical parameter $\nu$ of the anyons 
by $q = e^{i\pi\nu}$. The construction uses anyons contructed from usual 
fermionic oscillators and deformed bosonic oscillators. As a byproduct, 
realization deformed in any sector of the quantum superalgebras $\uqsa$ is 
obtained.
\end{abstract}

\vfill
\vfill

\rightline{CRM-2383}
\rightline{DSF-T-35/96}
\rightline{DFTT-39/96}
\rightline{ENSLAPP-AL-603/96}
\rightline{q-alg/9609033}
\rightline{September 1996}

\newpage
\pagestyle{plain}

\indent


\sect{Introduction}
\label{sectintro}

Superalgebras, which are the mathematical framework for describing symmetry
betwen bosons and fermions have by now found several interesting applications 
in physics, even if the {\em fundamental} supersymmetry between elementary
constituents of the matter has not yet been supported by experimental evidence.
A further enlarged concept of symmetry represented by quantum algebras has
shown up in a large number of areas in physics.The fusion of these new 
enlarged symmetry structures is a natural step and it leads to the so-called
$q$-superalgebras. 

Moreover, the connection between quantum algebras and generalized statistics 
has been pointed out in several contexts. {\em Anyons} are typical objects of 
generalized statistics whose importance in two-dimensional physics is 
relevant. They have been used to construct Schwinger-Jordan like realization 
of the deformed classical finite Lie algebras \cite{LS93,FMS94} and of the 
deformed affine Lie algebras of the unitary and symplectic series 
\cite{F3S96,Conf96}. So it is natural to ask which kind 
of oscillators are necessary to build realizations of $q$-superalgebras.
We will show that the deformation of the affine unitary superalgebras (and
therefore also of the  finite unitary superalgebras) can be realized by means
of anyons and of a new type of generalized statistics objects which satisfy 
braiding relations and which will be called {\em bosonic anyons} for reasons 
which will be clear from their definitions, see Sec. \ref{sect3}. 

Let us emphasize that the construction we propose may be interesting in the 
study of systems of correlated electrons. In fact the so-called $t-J$ model 
\cite{ZR88}, which has been suggested as an appropriate starting point for 
the theory of the high-temperature superconductivity, is supersymmetric for 
particular values of the coupling constants and of the chemical potential, 
the Hamiltonian commuting with a $su(1|2)$. Moreover in ref. \cite{LS93b} it 
has been shown that the $t-J$ model at the supersymmetric point can be 
written in terms of anyons, which gives a new realization of supersymmetry. 
Although in this model no deformation appears, it is conceivable that anyons 
can be used to describe further deformed generalizations of this or of 
similar models, like the Hubbard model \cite{Mon94}.

The article is organized in the following way:
in Sec. \ref{sect1} we briefly recall the structure of $\sa$ in the
Cartan-Weyl and in the Serre-Chevalley bases and then we write its deformation
in the distinguished Serre-Chevalley basis; in Sec. \ref{sect2} the 
fermionic-bosonic oscillators realization of $\sa$ is presented and finally 
in Sec. \ref{sect3} the anyonic realization of $\uqa$ is given in terms of 
one-dimensional anyons and bosonic anyons; in Sec. \ref{sect4} the 
generalization of the construction to two-dimensional anyons is discussed
and a few conclusions are presented.


\sect{Presentation of the superalgebra $\sa$} 
\label{sect1}

We will recall in this section the presentation of the affine Lie
superalgebra $\sa$, where $M,N \ge 1$, both in the Cartan--Weyl basis and in 
the Serre--Chevalley basis. We set $R=M+N-1$ and we exclude the case $R=1$
(obtained when $M=N=1$).

\subsection{Cartan--Weyl presentation of $\sa$}

In the Cartan--Weyl basis, the generators of the affine Lie superalgebra
$\sa$ are denoted by $h_a^m$ (Cartan generators) and $e_\fra^m$ (root
generators) where $a=1,\dots,R$ and $m \in \ZZ$. The even root system
of $\sa$ is given by $\Delta_0 = \{ \pm(\eps_i - \eps_j), \pm(\del_k
- \del_l) \}$  and the odd root system by $\Delta_1 = \{ \pm(\eps_i -
\del_k) \}$ where $1 \le i < j \le M$ and $1 \le k < l \le N$, the
$\eps_i$ and $\del_k$ spanning the dual of the Cartan subalgebra of
$gl(M|N)$. To each root generator $e_\fra^m$ one assigns a $\ZZ_2$-grading 
defined by $\deg(e_\fra^m)=0$ if $\fra \in \Delta_0$ and $\deg(e_\fra^m)=1$ 
if $\fra \in \Delta_1$. 
The generators satisfy for $M \ne N$ the following commutation relations:
\subequations
\bea
&& \cle h_a^m,h_b^n \cri = \gamma \,m \,\delta_{m+n,0}
\,K(h_a,h_b) \label{eq1a} \\
&& \cle h_a^m,e_\fra^n \cri = \fra_a ~ e_\fra^{m+n}
\label{eq1b} \\
&& \sle e_\fra^m,e_\frb^n \sri = 
\left\{ \begin{array}{ll}
\eps(\fra,\frb) ~ e_{\fra+\frb}^{m+n} & \quad \smbox{if} \fra + \frb
\smbox{is a root} \\ 
\fra^a h_a^{m+n} + \gamma \,m \,\delta_{m+n,0} \,K(e_\fra,e_{-\fra}) 
& \quad \smbox{if} \frb = -\fra \\ 
0 & \quad \smbox{otherwise}
\end{array} \right.
\label{eq1c} \\
&& \cle h_a^m,\gamma \cri = \cle e_\fra^m,\gamma \cri = 0
\label{eq1d}
\ena
\endsubequations 
where $\eps(\fra,\frb) = \pm 1$ is the usual 2-cocycle, $K$ is the 
(non-degenerate) Killing form on the horizontal superalgebra $A(M-1,N-1)$ 
and $\gamma$ is the central charge. $[\hspace{-2pt}[ ~,~ ]\hspace{-2pt}]$ 
denotes the super-commutator:
$[\hspace{-2pt}[ e_\fra^m,e_\frb^n ]\hspace{-2pt}] = e_\fra^m \, e_\frb^n 
- (-1)^{\deg(e_\fra^m).\deg(e_\frb^n)} e_\frb^n \, e_\fra^m$. Note that by 
virtue of Eqs. (\ref{eq1a}-d) the value of the central charge of 
$\hat A_{M-1}$ is opposite to that of $\hat A_{N-1}$.
\\
In the case $M=N$, although the Killing form is zero, it is possible to 
define a non-degenerate bilinear form $K$ on $A(N-1,N-1)$ such that Eqs. 
(\ref{eq1a}-d) still hold.

\subsection{Serre--Chevalley presentation of $\sa$}

In the Serre--Chevalley basis, the algebra is described in terms of
simple root and Cartan generators, the only data being the entries of the
Cartan matrix $(a_{\alpha\beta})$ of the algebra. Let us denote the
generators in the Serre--Chevalley basis by $h_\alpha$ and
$e_\alpha^\pm$ where $\alpha=0,1,\dots,R$. If $\tau$ is a subset of
$\{0,1,\dots,R\}$, the $\ZZ_2$-gradation of the superalgebra is defined
by setting $\deg(e_\alpha^\pm) = 0$ if $\alpha \notin \tau$ and 
$\deg(e_\alpha^\pm) = 1$ if $\alpha \in \tau$. The superalgebra is described 
by the (super)commutation relations
\subequations
\bea
&& \cle h_\alpha,h_\beta \cri = 0 \label{eq2a} \\
&& \cle h_\alpha,e^\pm_\beta \cri = \pm a_{\alpha\beta} e^\pm_\beta 
\label{eq2b} \\
&& \sle e^+_\alpha,e^-_\beta \sri = e^+_\alpha e^-_\beta 
- (-1)^{\deg(e^+_\alpha)\deg(e^-_\beta)} ~ e^-_\beta e^+_\alpha = 
h_\alpha ~ \delta_{\alpha\beta} \label{eq2c} \\
&& \ale e_\alpha^\pm,e_\alpha^\pm \ari = 0 
\medbox{if} a_{\alpha\alpha} = 0 \label{eq2d}
\ena
\endsubequations
and by the following relations:
\begin{itemize}
\item Serre relations for all $\alpha \ne \beta$
\be
(\mbox{ad}\, e_\alpha^\pm)^{1-{\tilde a}_{\alpha\beta}} ~ e_\beta^\pm = 0
\label{eq3}
\ee
\item supplementary relations for $\alpha$ such that $a_{\alpha\alpha}=0$
\be
\sle (\mbox{ad}\,e_{\alpha-1}^\pm) \,e_\alpha^\pm, 
(\mbox{ad}\,e_{\alpha+1}^\pm) \,e_\alpha^\pm \sri = 0
\label{eq4}
\ee
\end{itemize}
where the matrix ${\tilde A} = ({\tilde a}_{\alpha\beta})$ is deduced
{from} the Cartan matrix $A = (a_{\alpha\beta})$ of $\sa$ by replacing
all its positive off-diagonal entries by $-1$.
Here ad denotes the adjoint action: 
$(\mbox{ad}\,X) ~ Y = XY - (-1)^{\deg X.\deg Y} ~ YX$.

\medskip

One has to emphasize that for superalgebras, the description given by the 
Serre relations (\ref{eq3}) leads in general to a bigger superalgebra than
the superalgebra under consideration. It is necessary to write supplementary 
relations involving more than two generators, that for $A(M-1,N-1)$ take the 
form (\ref{eq4}), in order to quotient the bigger superalgebra and recover 
the original one, see \cite{Sch93} for more details. 
As one can imagine, these supplementary conditions appear when one deals with 
isotropic fermionic simple roots (that is $a_{\alpha\alpha} = 0$). Note that 
these supplementary relations are unnecessary when $M=1$ or $N=1$.
 
\medskip

In the following, we will only use the Serre--Chevalley description of the 
affine Lie superalgebra in the {\em distinguished} basis, such that the 
number of odd simple roots is the smallest one. In the case of $\sa$, the
distinguished basis is defined by taking $\tau = \{0,M\}$. The
corresponding Dynkin diagram is the following, with the labels identifying 
the corresponding simple roots:
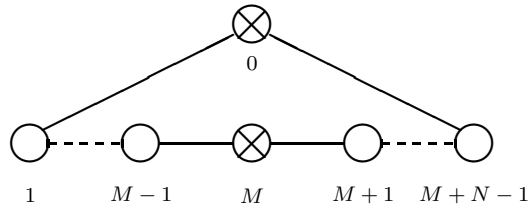
\begin{center}
\begin{picture}(180,80)
\thicklines
\multiput(0,20)(42,0){5}{\circle{14}}
\put(84,65){\circle{14}}
\put(7,20){\dashbox{3}(28,0)}
\put(49,20){\line(1,0){28}}
\put(91,20){\line(1,0){28}}
\put(5,25){\line(2,1){72}}
\put(163,25){\line(-2,1){72}}
\put(133,20){\dashbox{3}(28,0)}
\put(79,15){\line(1,1){10}}\put(79,25){\line(1,-1){10}}
\put(79,60){\line(1,1){10}}\put(79,70){\line(1,-1){10}}
\put(0,0){\makebox(0.4,0.6){\scriptsize{$1$}}}
\put(42,0){\makebox(0.4,0.6){\scriptsize{$M-1$}}}
\put(84,0){\makebox(0.4,0.6){\scriptsize{$M$}}}
\put(126,0){\makebox(0.4,0.6){\scriptsize{$M+1$}}}
\put(168,0){\makebox(0.4,0.6){\scriptsize{$M+N-1$}}}
\put(84,50){\makebox(0.4,0.6){\scriptsize{$0$}}}
\end{picture}
\end{center}
associated to the Cartan matrix
\be
\left(\begin{array}{rrrrrrrrrrrr}
0 & 1 & 0 & \cdots & \cdots &&&& \cdots & \cdots & 0 & -1 \\
-1 & 2 & -1 & 0 &&&&&&&& 0 \\
0 & -1 & \ddots & \ddots & \ddots &&&&&&& \vdots \\
\vdots & 0 & \ddots && \ddots & 0 &&&&&& \vdots \\
\vdots  && \ddots & \ddots & \ddots & -1 & \ddots &&&&& \\
&&& 0 & -1 & 2 & -1 & \ddots &&&& \\ 
&&&& \ddots & -1 & 0 & 1 & \ddots &&& \\ 
&&&&& \ddots & -1 & 2 & -1 & 0 && \vdots \\
\vdots &&&&&& \ddots & -1 & \ddots & \ddots & \ddots & \vdots \\
\vdots &&&&&&& 0 & \ddots && \ddots & 0 \\
0 &&&&&&&& \ddots & \ddots & \ddots & -1 \\
-1 & 0 & \cdots & \cdots && \cdots & \cdots & 0 & \cdots & 0 & -1 & 2 \\
\end{array}\right)
\label{eq5}
\ee
The correspondence between the distinguished Serre--Chevalley and the
Cartan--Weyl bases is the following ($i=1,\dots,M-1$ and 
$k=1,\dots,N-1$): 
\be
\begin{array}{lll}
\bigg. h_i = h_i^0 \qquad & 
e_i^+ = e_{\eps_i-\eps_{i+1}}^0 \qquad & 
e_i^- = e_{\eps_{i+1}-\eps_i}^0 \\
\bigg. h_M = h_M^0 \qquad & 
e_M^+ = e_{\eps_M-\del_1}^0 \qquad & 
e_M^- = e_{\del_1-\eps_M}^0 \\
\bigg. h_{M+k} = h_{M+k}^0 \qquad & 
e_{M+k}^+ = e_{\del_k-\del_{k+1}}^0 \qquad & 
e_{M+k}^- = e_{\del_{k+1}-\del_k}^0 \\
\bigg. h_0 = -\gamma + \sum_{i=1}^{M} h_i^0 - \sum_{k=1}^{N-1} h_{M+k}^0
\qquad 
& e_0^+ = e_{\del_N-\eps_1}^1 \qquad & 
e_0^- = e_{\eps_1-\del_N}^{-1}
\end{array}
\label{eq6}
\ee
Notice that in the Serre--Chevalley picture the central charge $\gamma$ \
is uniquely defined by the following equation:
\be
h_0 = -\gamma + \sum_{i=1}^{M} h_i - \sum_{k=1}^{N-1} h_{M+k}
\label{eq6bis}
\ee

\subsection{Serre--Chevalley presentation of $\uqa$}

We consider now the universal quantum affine Lie superalgebra $\uqa$. The 
Serre--Chevalley description in the quantum case is very similar. The 
defining relations take the form
\subequations
\bea
&& \cle h_\alpha,h_\beta \cri = 0 \label{eq7a} \\
&& \cle h_\alpha,e_\beta^\pm \cri = \pm a_{\alpha\beta} e_\beta^\pm 
\label{eq7b} \\ 
&& \sle e_\alpha^+,e_\beta^- \sri = e^+_\alpha e^-_\beta 
- (-1)^{\deg(e^+_\alpha)\deg(e^-_\beta)} ~ e^-_\beta e^+_\alpha =
\frac{q_\alpha^{h_\alpha} - q_\alpha^{-h_\alpha}}{q_\alpha - q_\alpha^{-1}} ~
\delta_{\alpha\beta} \label{eq7c} \\ 
&& \ale e_\alpha^\pm,e_\alpha^\pm \ari = 0 
\medbox{if} a_{\alpha\alpha} = 0 \label{eq7d}
\ena
\endsubequations
where $q_{\alpha} = q^{d_\alpha}$ and the numbers $d_\alpha$ symmetrize the 
Cartan matrix $\bar a_{\alpha\beta}$ of $A(M-1,N-1)$: $d_\alpha \bar 
a_{\alpha\beta} = d_\beta \bar a_{\beta\alpha}$ ($\alpha,\beta \ne 0$) and 
$d_0 = 1$ (in the distinguished basis, $q_\alpha = q$ for $\alpha=0,\dots,M$ 
and $q_\alpha = q^{-1}$ for $\alpha=M+1,\dots,M+N-1$).
\\
In terms of the generators 
$\cE_\alpha^\pm = e_\alpha^\pm ~ q_\alpha^{-h_\alpha/2}$
the usual Serre relations are for all $\alpha \ne \beta$
\be
(\mbox{ad}_q\,\cE_\alpha^\pm)^{1-{\tilde a}_{\alpha\beta}} ~ \cE_\beta^\pm = 0
\label{eq8}
\ee
while the supplementary relations read now  for $\alpha$ such that 
$a_{\alpha\alpha}=0$ (the definition of the quantum adjoint action 
$\mbox{ad}_q$ is given below Eq. \ref{eq12}) \cite{FLV91,Sch93,Yam94}
\be 
\sle (\mbox{ad}_q\,\cE_{\alpha-1}^\pm) \,\cE_\alpha^\pm, 
(\mbox{ad}_q\,\cE_{\alpha+1}^\pm) \,\cE_\alpha^\pm \sri = 0
\label{eq10}
\ee
or in terms of the generators $e_\alpha^\pm$
\be
\sle \cle e_{\alpha-1}^\pm , e_\alpha^\pm \cri_q , 
\cle e_\alpha^\pm , e_{\alpha+1}^\pm \cri_q \sri = 0
\label{eq9}
\ee
the $q$-commutator being defined as usual by $\cle X,Y \cri_q = XY - q YX$.

\medskip

The universal quantum affine Lie superalgebra $\cU \equiv \uqa$ is endowed 
with a Hopf algebra structure, with coproduct $\Delta: \cU \rightarrow \cU 
\otimes \cU$, counit $\eps: \cU \rightarrow \CC$ and antipode $S: \cU 
\rightarrow \cU$ such that ($\alpha=0,1,\dots,R$)
\subequations
\bea
&& \Delta(h_\alpha) = 1 \otimes h_\alpha + h_\alpha \otimes 1 \medbox{and}
\Delta(e_\alpha^\pm) = e_\alpha^\pm \otimes q_\alpha^{h_\alpha/2} 
+ q_\alpha^{-h_\alpha/2} \otimes e_\alpha^\pm
\label{eq11a} \\
&& \eps(h_\alpha) = \eps(e_\alpha^\pm) = 0 \medbox{and} \eps(1) = 1
\label{eq11b} \\
&& S(h_\alpha) = -h_\alpha \medbox{and} S(e_\alpha^\pm) = 
-q_\alpha^{\pm a_{\alpha\alpha}/2} e_\alpha^\pm
\label{eq11c}
\ena
\endsubequations
The quantum adjoint action $\mbox{ad}_q$ can be explicitly written in terms 
of the coproduct and the antipode as
\be
(\mbox{ad}_q\,X) ~ Y = (-1)^{\deg X_{(2)}.\deg Y} ~ X_{(1)}\,Y\,S(X_{(2)})
\label{eq12}
\ee
using the Sweedler notation for the coproduct: 
$\Delta(X) = X_{(1)} \otimes X_{(2)}$ (summation is understood).


\sect{Oscillator realization of the affine Lie superalgebra $\sa$}
\label{sect2}

Let us recall now the oscillator realization of $\sa$ in terms of
creation and annihilation operators. We consider an infinite number of
fermionic oscillators $c_i(r), c_i^\dagger(r)$ with $i=1,\dots,M$ and $r
\in \ZZ+1/2 = \ZP$, which satisfy the anticommutation relations
\be
\ale c_i(r),c_j(s) \ari = \ale c_i^\dagger(r), c_j^\dagger(s) \ari = 0 
\bigbox{and} 
\ale c_i(r),c_j^\dagger(s) \ari = \delta_{ij} \delta_{rs} 
\label{eq20}
\ee
and an infinite number of bosonic oscillators $d_k(r), d_k^\dagger(r)$
with $k=1,\dots,N$ and $r \in \ZP$, which satisfy the commutation 
relations
\be
\cle d_k(r),d_l(s) \cri = \cle d_k^\dagger(r), d_l^\dagger(s) \cri = 0 
\bigbox{and} 
\cle d_k(r),d_l^\dagger(s) \cri = \delta_{kl} \delta_{rs} 
\label{eq21}
\ee
the two sets $c_i(r), c_i^\dagger(r)$ and $d_k(r), d_k^\dagger(r)$
commuting each other:
\be
\cle c_i(r),d_k(s) \cri = \cle c_i(r),d_k^\dagger(s) \cri = 
\cle c_i^\dagger(r),d_k(s) \cri = \cle c_i^\dagger(r),d_k^\dagger(s) 
\cri = 0 \label{eq30}
\ee
The fermionic and bosonic number operators are defined as usual by 
$n_i(r) = c_i^\dagger(r)c_i(r)$ and $n'_k(r) = d_k^\dagger(r)d_k(r)$. 
\\
These oscillators are equipped with a normal ordered product such that
\be
:c_i^\dagger(r) c_j(s): = \left\{ \begin{array}{ll}
c_i^\dagger(r) c_j(s) & \smbox{if} s > 0 \\
- c_j(s) c_i^\dagger(r) & \smbox{if} s < 0
\end{array} \right.
\label{eq22}
\ee
and
\be
:d_k^\dagger(r) d_l(s): = \left\{ \begin{array}{ll}
d_k^\dagger(r) d_l(s) & \smbox{if} s > 0 \\
d_l(s) d_k^\dagger(r) & \smbox{if} s < 0
\end{array} \right.
\label{eq23}
\ee
Therefore
\be
:n_i(r): = \left\{ \begin{array}{ll}
n_i(r)  & \smbox{if} r > 0 \\
n_i(r) - 1  & \smbox{if} r < 0
\end{array} \right.
\label{eq24}
\ee
and
\be
:n'_k(r): = \left\{ \begin{array}{ll}
n'_k(r)  & \smbox{if} r > 0 \\
n'_k(r) + 1  & \smbox{if} r < 0 \\
\end{array} \right.
\label{eq25}
\ee
Then an oscillator realization of the generators of $\sa$ in the
Cartan--Weyl basis with $\gamma =1$ is given by 
\subequations
\bea
&& \hspace{-5mm} h_i^m = \sum_{r \in \ZP} \left(:c_i^\dagger(r) c_i(r+m): - 
:c_{i+1}^\dagger(r) c_{i+1}(r+m): \right) ~~~ (i=1,\dots,M-1) \,, 
\label{eq26a} \\ 
&& \hspace{-5mm} h_M^m = \sum_{r \in \ZP} \left(:c_M^\dagger(r) c_M(r+m): + 
:d_1^\dagger(r) d_1(r+m): \right) \,, \label{eq26b} \\ 
&& \hspace{-5mm} h_{M+k}^m = \sum_{r \in \ZP} \left(:d_k^\dagger(r) d_k(r+m): \- :d_{k+1}^\dagger(r) d_{k+1}(r+m): \right) ~~~ (k=1,\dots,N-1) \,, 
\label{eq26c} \\ 
&& \hspace{-5mm} e_{\eps_i-\eps_j}^m = \sum_{r \in \ZP} 
c_i^\dagger(r) c_j(r+m) \label{eq26d} \\
&& \hspace{-5mm} e_{\del_k-\del_l}^m = \sum_{r \in \ZP} 
d_k^\dagger(r) d_l(r+m) \label{eq26e} \\
&& \hspace{-5mm} e_{\eps_i-\del_k}^m = \sum_{r \in \ZP} 
c_i^\dagger(r) d_k(r+m) \label{eq26f} \\
&& \hspace{-5mm} e_{\del_k-\eps_i}^m = \sum_{r \in \ZP} 
d_k^\dagger(r) c_i(r+m) \label{eq26g} 
\ena
\endsubequations
A fermionic oscillator realization of the simple generators of $\sa$ in
the distinguished Serre-Chevalley basis is given by
$(\alpha=0,1,\dots,R)$
\bea
h_\alpha = \sum_{r \in \ZP} h_\alpha(r)
\bigbox{and}
e_\alpha^\pm = \sum_{r \in \ZP} e_\alpha^\pm (r)
\label{eq27} 
\ena
where ($i=1,\dots,M-1$ and $k=1,\dots,N-1$)
\subequations
\bea
&& h_i(r) = n_i(r) - n_{i+1}(r) = ~ :n_i(r): - :n_{i+1}(r):
\label{eq28a} \\
&& h_M(r) = n_M(r) + n'_1(r) = ~ :n_M(r): + :n'_1(r):
\label{eq28b} \\
&& h_{M+k}(r) = n'_k(r) - n'_{k+1}(r) = ~ :n'_k(r): - :n'_{k+1}(r):
\label{eq28c} \\
&& h_0(r) = n'_N(r) + n_1(r+1) = ~ :n'_N(r): + :n_1(r+1): -
~ \delta_{r+1/2,0} \label{eq28d} \\
&& e_i^+(r) =  c_i^\dagger(r) c_{i+1}(r) \,, \hspace{18mm} 
e_i^-(r) =  c_{i+1}^\dagger(r) c_i(r) \label{eq28e} \\
&& e_M^+(r) =  c_M^\dagger(r) d_1(r) \,, \hspace{18.5mm} 
e_M^-(r) =  d_1^\dagger(r) c_M(r) \label{eq28f} \\
&& e_{M+k}^+(r) =  d_k^\dagger(r) d_{k+1}(r) \,, \hspace{12mm} 
e_{M+k}^-(r) =  d_{k+1}^\dagger(r) d_k(r) \label{eq28g} \\
&& e_0^+(r) =  d_N^\dagger(r) c_1(r+1) \,, \hspace{14mm}
e_0^-(r) =  c_1^\dagger(r+1) d_N(r) \label{eq28h} 
\ena
\endsubequations
Inserting Eq. (\ref{eq28d}) into Eq. (\ref{eq27}) and taking into
account that the sum over $r$ can be splitted into a sum of two
convergent series only after normal ordering, one can check that 
\be
h_0 = -1 + \sum_{r \in \ZP} :n'_N(r): + \sum_{r \in \ZP} :n_1(r): ~ 
= -1 + \sum_{i=1}^{M-1} h_i^0 + h_M^0 - \sum_{k=1}^{N-1} h_{M+k}^0 \,,
\label{eq29}
\ee
that is the central charge is $\gamma = 1$.
\\
Note that the value of the central charge is related to the definition
of the normal ordered product. A different definition like ($i = 1,\dots,M$
and $k=1,\dots,N$)
\be
:n_i(r): = n_i(r) \smbox{and} :n'_k(r): = n'_k(r) 
\smbox{for any} r \in \ZP
\label{eq31}
\ee
would lead to $\gamma = 0$.


\sect{Anyonic realization of $\uqa$}
\label{sect3}

In order to obtain an anyonic realization of $\uqa$, we will replace the
fermionic and bosonic oscillators by suitable anyons in the expressions
of the simple generators of $\uqa$ in the distinguished Serre--Chevalley
basis. Since we have to deal with fermionic and bosonic generators, we
have to introduce two different types of anyons.

\medskip

Let us first define fermionic anyons on a one-dimensional lattice $\ZP$ 
\cite{LS93,FLS96}:
\be
a_i(r) = K_i(r) c_i(r) \bigbox{and} \ta_i(r) = \tK_i(r) c_i(r)
\label{eq40}
\ee
and similar expressions for the conjugated operators $a_i^\dagger(r)$ 
and $\ta_i^\dagger(r)$, where the disorder factors $K_i(r)$ and
$\tK_i(r)$ are expressed as
\subequations
\bea
&& K_i(r) = q^{-\half\sum_{t\in\ZP} \eps(t-r) :n_i(t):} \label{eq41a} \\
&& \tK_i(r) = q^{\half\sum_{t\in\ZP} \eps(t-r) :n_i(t):} \label{eq41b}
\ena
\endsubequations
The function $\eps(t) = |t|/t$ if $t \ne 0$ and $\eps(0) = 0$ is the 
sign function.
\\
It is easy to prove that the $a$-type anyons satisfy the following
braiding relations for $r>s$:
\bea
&& a_i(r) a_i(s) + q^{-1} a_i(s) a_i(r) = 0 \nonumber \\
&& a_i^\dagger(r) a_i^\dagger(s) + q^{-1} a_i^\dagger(s)
a_i^\dagger(r) = 0 \nonumber \\ 
&& a^\dagger_i(r) a_i(s) + q ~ a_i(s) a^\dagger_i(r) = 0 \nonumber \\
&& a_i(r) a_i^\dagger(s) + q ~ a_i^\dagger(s) a_i(r) = 0 
\label{eq42}
\ena
and
\bea
&& a_i(r) a_i^\dagger(r) + a_i^\dagger(r) a_i(r) = 1 \nonumber\\
&& a_i(r)^2 =  a_i^\dagger(r)^2 = 0 
\label{eq43}
\ena
The braiding relations between the $\ta$-type anyons are obtained from
eqs. (\ref{eq42}) and (\ref{eq43}) by replacing $q \leftrightarrow
q^{-1}$.
\\
Finally, the braiding relations between $a$-type and $\ta$-type anyons
are given by
\bea
&& \ale {\ta}_i(r), a_i(s) \ari = \ale {\ta}^\dagger_i(r),a^\dagger_i(s) 
\ari = 0 \medbox{for all} r,s \in \ZP 
\label{eq44} \\
&& \ale {\ta}^\dagger_i(r), a_i(s) \ari = \ale {\ta}_i(r),a^\dagger_i(s) 
\ari = 0 \medbox{for all} r \ne s \in \ZP 
\label{eq45}
\ena
and
\bea
&& \ale {\ta}_i(r), a^\dagger_i(r) \ari =
q^{ ~ \sum_{t \in \ZP} \eps(t-r) :n_i(t): } \nonumber \\
&& \ale {\ta}^\dagger_i(r), a_i(r) \ari = 
q^{ - \sum_{t \in \ZP} \eps(t-r) :n_i(t): } 
\label{eq46}
\ena
Moreover, the following identity holds:
\be
a^\dagger_i(r) a_i(r) =  {\ta}^\dagger_i(r) {\ta}_i(r) = n_i(r)
\label{eq47}
\ee
the normal ordering between $a$-type and $\ta$-type anyons being defined
as in eq. (\ref{eq22}).

\medskip

Now we will define anyonic-like operators based on $q$-deformed bosons.
Let us recall that $q$-deformed bosons can be constructed from ordinary
ones by the following procedure \cite{Son90}:
\subequations
\bea
&& n'_k(r) = d_k^\dagger(r) d_k(r) \label{eq48a} \\
&& b_k(r) = d_k(r) ~ \sqrt{\frac{[n'_k(r)]_q}{n'_k(r)}} = 
\sqrt{\frac{[n'_k(r)+1]_q}{n'_k(r)+1}} ~ d_k(r) \label{eq48b} \\
&& b_k^\dagger(r) = \sqrt{\frac{[n'_k(r)]_q}{n'_k(r)}} ~ d_k^\dagger(r) = 
d_k^\dagger(r) ~ \sqrt{\frac{[n'_k(r)+1]_q}{n'_k(r)+1}} \label{eq48c}
\ena
\endsubequations
The $q$-deformed bosons $b_k(r),b_k^\dagger(r)$ satisfy the following
$q$-commutation relations: 
\subequations
\bea
&& b_k(r) b_l^\dagger(s) - q^{\delta_{kl}\delta_{rs}} b_l^\dagger(s)
b_k(r) =  q^{-n'_k(r)} \delta_{kl} \delta_{rs} \label{eq49a} \\
&& b_k(r) b_l^\dagger(s) - q^{-\delta_{kl}\delta_{rs}} b_l^\dagger(s)
b_k(r) =  q^{n'_k(r)} \delta_{kl} \delta_{rs} \label{eq49b} \\
&& b_k(r) b_l(s) - b_l(s) b_k(r) = b_k^\dagger(r) b_l^\dagger(s) 
- b_l^\dagger(s) b_k^\dagger(r) = 0 \label{eq49c} \\
&& \cle n'_k(r),b_l(s) \cri = -b_k(r) \delta_{kl}\delta_{rs} \label{eq49d} \\ 
&& \cle n'_k(r),b_l^\dagger(s) \cri = b_k^\dagger(r) \delta_{kl}\delta_{rs} 
\label{eq49e}
\ena
\endsubequations
{from} which it follows that
\be
b_k^\dagger(r) b_k(r) = [n'_k(r)]_q \bigbox{and}
b_k(r) b_k^\dagger(r) = [n'_k(r)+1]_q \label{eq50}
\ee
Now, let us define anyonic-like operators as follows:
\be
A_k(r) = K'_k(r) b_k(r) \bigbox{and} \tA_k(r) = \tK'_k(r) b_k(r)
\label{eq51}
\ee
and similar expressions for the conjugated operators $A_k^\dagger(r)$ 
and $\tA_k^\dagger(r)$, where the disorder factors are given by
\subequations
\bea
&& K'_k(r) = q^{\half\sum_{t\in\ZP} \eps(t-r) :n'_k(t):} \label{eq52a} \\
&& \tK'_k(r) = q^{-\half\sum_{t\in\ZP} \eps(t-r) :n'_k(t):} \label{eq52b}
\ena
\endsubequations
It can be proved that the operators $A_k(r),A_k^\dagger(r)$ satisfy the
following braiding relations for $r>s$:
\bea
&& A_k(r) A_k(s) - q A_k(s) A_k(r) = 0 \nonumber \\
&& A_k^\dagger(r) A_k^\dagger(s) - q A_k^\dagger(s)
A_k^\dagger(r) = 0 \nonumber \\ 
&& A^\dagger_i(r) A_k(s) - q^{-1} ~ A_k(s) A^\dagger_i(r) = 0 
\nonumber \\ 
&& A_k(r) A_k^\dagger(s) - q^{-1} ~ A_k^\dagger(s) A_k(r) = 0 
\label{eq53}
\ena
and
\bea
&& A_k(r) A_k^\dagger(r) - q A_k^\dagger(r) A_k(r) = q^{-n'_k(r)}
\nonumber \\
&& A_k(r) A_k^\dagger(r) - q^{-1} A_k^\dagger(r) A_k(r) = q^{n'_k(r)}
\label{eq54}
\ena
Therefore, the operators $A_k(r),A_k^\dagger(r)$ satisfy the
$q$-commutation relations of the $q$-deformed bosonic oscillator at the
same point, while they satisfy braiding relations when taken at different
points.
\\
The braiding relations between the $\tA$-type anyons are obtained from
eq. (\ref{eq53}) by replacing $q \leftrightarrow q^{-1}$.

Note that for these anyonic $q$-deformed bosons defined by the above 
relations, we do not have any physical interpretations, on the contrary of 
the $a$-type anyons, see refs. \cite{LS93,FLS96}.
Is is worth to point out that the above introduced {\em bosonic anyons}
differ from the ones introduced in ref. \cite{LMR96} by the non trivial fact
that our anyons are defined on a lattice while in ref. \cite{LMR96} are defined
in the continuum and by the local braiding relation (\ref{eq54}). Replacing 
in Eq. (\ref{eq51}) the $q$-boson by a standard boson, we find in the lattice 
the bosonic anyons of ref. \cite{LMR96}. We will come back on the difference 
between the two approaches in the next section.

\medskip

Now we can build an anyonic realization of $\uqa$ by "anyonizing" the
oscillator realization eqs. (\ref{eq28a}-\ref{eq28h}), that is
replacing the fermionic oscillators $c_i$ by the anyonic oscillators
$a_i$ and $\ta_i$ and the bosonic oscillators $b_i$ by the operators
$A_i$ and $\tA_i$. More precisely, one has:

\begin{theorem}
An anyonic realization of the simple generators of the quantum affine Lie
superalgebra $\uqa$ with central charge $\gamma = 1$ is given by (with
$\alpha=0,1,\dots,R$) 
\be
H_\alpha = \sum_{r \in \ZP} H_\alpha(r)
\bigbox{and} 
E_\alpha^\pm = \sum_{r \in \ZP} E_\alpha^\pm (r)
\label{eq55} 
\ee
where ($i=1,\dots,M-1$ and $k=1,\dots,N-1$)
\subequations
\bea
&& H_i(r) = n_i(r) - n_{i+1}(r) =  :n_i(r): - :n_{i+1}(r):
\label{eq56a} \\
&& H_M(r) = n_M(r) + n'_1(r) =  :n_M(r): + :n'_1(r):
\label{eq56b} \\
&& H_{M+k}(r) = n'_k(r) - n'_{k+1}(r) =  :n'_k(r): - :n'_{k+1}(r):
\label{eq56c} \\
&& H_0(r) = n'_N(r) + n_1(r+1) =  :n'_N(r): + :n_1(r+1): -
\delta_{r+1/2,0} \label{eq56d} \\
&& E_i^+(r) =  a_i^\dagger(r) a_{i+1}(r) \,, \hspace{20mm} 
E_i^-(r) =  \ta_{i+1}^\dagger(r) \ta_i(r) \label{eq56e} \\
&& E_M^+(r) =  a_M^\dagger(r) A_1(r) \,, \hspace{19mm} 
E_M^-(r) =  \tA_1^\dagger(r) \ta_M(r) \label{eq56f} \\
&& E_{M+k}^+(r) =  A_k^\dagger(r) A_{k+1}(r) \,, \hspace{12mm} 
E_{M+k}^-(r) =  \tA_{k+1}^\dagger(r) \tA_k(r) \label{eq56g} \\
&& E_0^+(r) =  A_N^\dagger(r) a_1(r+1) \,, \hspace{14mm}
E_0^-(r) =  \ta_1^\dagger(r+1) \tA_N(r) \label{eq56h} 
\ena
\endsubequations
\end{theorem}

\noindent
{\bf Proof}
We must check that the realization eqs. (\ref{eq55}) and
(\ref{eq56a}-\ref{eq56h}) indeed satisfy the quantum affine Lie
superalgebra $\uqa$ in the distinguished Serre--Chevalley basis
(\ref{eq7a}-\ref{eq7d}) together with the quantum Serre relations
(\ref{eq8}) and (\ref{eq9}). The proof follows the lines of the algebraic
case \cite{F3S96}: the equations (\ref{eq7a}-\ref{eq9}) which define a
generic deformed affine superalgebra $\cU_q(\widehat\cA)$ reduce to 
$\cU_q(\cA)$ when the affine dot is removed and to another finite
dimensional superalgebra $\cU_q(\cA')$ if the affine dot is kept and
one or more other suitable dots are removed. The relations defining
$\cU_q(\widehat\cA)$ coincide with the union of those defining $\cU_q(\cA)$
and $\cU_q(\cA')$: therefore, it will be enough to check that the equations
defining $\cU_q(\cA)$ and $\cU_q(\cA')$ are satisfied.

Consider the non-extended Dynkin diagram of $A(M-1,N-1)$ to which
the set of generators $\{ H_\alpha, E^\pm_\alpha \}$ (with $\alpha \ne 0$) 
corresponds. Inserting Eqs. (\ref{eq40}), (\ref{eq41a}-b), (\ref{eq50}),
(\ref{eq51}), the expressions (\ref{eq56e}-g) become
\be
E^\pm_\alpha(r) = \hat e^\pm_\alpha(r) ~ q_\alpha^{\half \sum_{t \in \ZP} 
\eps(t-r) \, :h_\alpha(t):}
\label{eq57}
\ee
where the generators $\hat e^\pm_\alpha(r)$ are obtained from the generators 
$e^\pm_\alpha(r)$ in Eqs. (\ref{eq28e}-g) replacing the bosonic oscillators 
$d_k$ by the $q$-deformed bosons $b_k$. The generators $\hat e^\pm_\alpha(r)$ 
coincide locally, i.e. for fixed $r$, with the generators of $\uqsa$ of refs. 
\cite{FLV91,FSV91} as the $q$-deformed fermions $\psi_i$ of \cite{FSV91} are 
equivalent to the usual fermionic oscillators $c_i$.
It follows that the generators $\{ H_\alpha, E^\pm_\alpha \}$ of 
Eq. (\ref{eq55}) are a representation of $\uqsa$ as they are obtained with 
the correct coproduct, see Eqs. (\ref{eq11a}) and (\ref{eq57}), 
by the representation in terms of $\{ h_\alpha, \hat e^\pm_\alpha \}$.
Let us remark that due to the equivalence betwen $q$-fermions and standard 
fermions the realization of finite $q$-superalgebras of ref. \cite{FSV91} 
are realizations of deformed algebras only for the subalgebra realized 
in terms of $q$-bosons while the subalgebra realized in terms of $q$-fermions 
is left undeformed. On the contrary the here presented anyonic realization 
is completely deformed in any sector. 
Finally let us remark that the difference $q \rightarrow q^{-1}$ in the 
disorder factor of the $A$-type anyons (in the site $r$) with respect to the 
$a$-type anyons (in the same site) -- see the $q_\alpha$-factor in 
Eq. (\ref{eq57}) -- is essential for the consistency of the 
$q$-superalgebra structure.

We consider then the {\em extended} Dynkin diagram of $A(M-1,N-1)$ and
we delete a dot which is not the affine dot. For example, cutting the
dot number 2, we obtain the following Dynkin diagram:
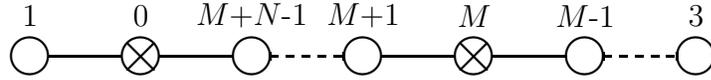
\begin{center}
\begin{picture}(280,20)
\thicklines
\multiput(0,0)(42,0){7}{\circle{14}}
\put(0,15){\makebox(0.4,0.6){1}}
\put(42,15){\makebox(0.4,0.6){0}}
\put(84,15){\makebox(0.4,0.6){$M$+$N$-1}}
\put(126,15){\makebox(0.4,0.6){$M$+1}}
\put(168,15){\makebox(0.4,0.6){$M$}}
\put(210,15){\makebox(0.4,0.6){$M$-1}}
\put(252,15){\makebox(0.4,0.6){3}}
\put(37,-5){\line(1,1){10}}\put(37,5){\line(1,-1){10}}
\put(163,-5){\line(1,1){10}}\put(163,5){\line(1,-1){10}}
\put(7,0){\line(1,0){28}}
\put(49,0){\line(1,0){28}}
\put(91,0){\dashbox{3}(28,0)}
\put(133,0){\line(1,0){28}}
\put(175,0){\line(1,0){28}}
\put(217,0){\dashbox{3}(28,0)}
\end{picture}
\end{center}
which corresponds to the Lie superalgebra $A(M-1,N-1)$ in a {\em
non-distinguished} basis.
\\
For a fixed $r \in \ZP$, it is possible to show that the set 
$\{ h_j(r), \hat e^\pm_j(r) , h_1(r+1), \hat e^\pm_1(r+1) \}$ 
($j = 0,3,\dots,M+N-1$) is a representation of 
$\uqsa$ in the non-distinguished basis specified by the above Dynkin diagram. 
We emphasize that in this case we have to satisfy two more supplementary Serre
relations than in the distinguished basis. Of course for particular values
of $M$ and $N$ one or both relations can be absent. Note that deleting the
$M$-th dot, we reobtain the superalgebra $\uqsa$ in the distinguished
basis again. Then it follows that the
generators $\{ H_\alpha, E^\pm_\alpha \}$ with $\alpha \ne 2$ are
a representation of $\uqsa$ as they are obtained by the generators of
a representation of the finite $q$-superalgebra with the correct coproduct.
This completes the proof.


\sect{General representations and conclusions}
\label{sect4}

In the previous section, we have built a representation of the
deformed affine Lie superalgebras $\uqa$ by means of anyons defined on an
infinite linear chain; as the corresponding fermionic representation,
it has central charge $\gamma = 1$.
Representations with vanishing central charge could be built in the
same way by using alternative normal ordering prescriptions Eq. (\ref{eq31}).

Representations with $\gamma=0$ and $\gamma=1$ can be combined together to 
get representations with arbitrary positive integer central charges 
(we do not discuss here the problem of the irreducibility of
these representations). Associating a representation to any horizontal line 
of a two-dimensional square lattice, infinite in one direction (say the 
horizontal one), and taking $K$ copies of representations in one-dimensional 
lattice with central charge equal to 1, one can get representations with the 
value of the central charge equal to $K$. Note that by combining one 
representation with central charge equal to $K$ with a finite number of 
representations (in one-dimensional lattice) with vanishing value of the 
central charge one obtains an inequivalent representation with the same value 
of the central charge. The extensions to two-dimensional lattice infinite in 
both directions can also be done, but it requires some care in the definition 
in order to avoid convergence problems.

\medskip

We have shown in ref. \cite{F3S96} that the use of $a$-anyons on a 
two-dimensional lattice naturally gives the correct coproduct with the 
correct powers of the deformation of the representations of a $q$-algebra 
defined in a fixed site of the lattice. For completeness we recall that each 
site of the two-dimensional lattice is labelled by a vector $\vec{x} = 
(x_1,x_2)$, the first component $x_1 \in \ZP$ being the coordinate of a site 
on the line $x_2 \in \ZZ$. The angle $\Theta(\vec{x},\vec{y})$ which enters 
in the definition of two-dimensional $a$-anyons through the disorder factor, 
see e.g. ref. \cite{FLS96}, 
\be
K(\vec{x}) = \exp\left(i\nu{\sum\limits_{\vec{y}\ne \vec{x}} 
\Theta(\vec{x},\vec{y}) ~ n(\vec{y})} \right)
\label{eq61}
\ee
may be chosen in such a way that
\be
\Theta(\vec{x},\vec{y}) = \left\{ \begin{array}{ll} 
+\pi/2 & \qquad \smbox{if} x_2 > y_2 \cr
-\pi/2 & \qquad \smbox{if} x_2 < y_2 
\end{array} \right.
\label{eq60}
\ee
while if $\vec{x}$ and $\vec{y}$ lie on the same horizontal line, that is
$x_2 = y_2$, the definitions of Sec. \ref{sect3} hold.
Two-dimensional anyons still satisfy the braiding and anticommutations 
relations expressed in the general form in Eqs. (\ref{eq42})-(\ref{eq46}).
Analogous relations hold for the $A$-anyons.

Let us replace in the equations of Sec. \ref{sect3} the one-dimensional 
anyons by two-dimensional ones and sum over the sites of the two-dimensional 
lattice. This sum has to be read as a sum over the infinite line $x_1$
and a sum over the finite number of lines labelled by $x_2$. Then the 
generators are given by a sum, with the correct coproduct, of the generators
of a $\uqa$ representations defined in a line. Therefore they define a $\uqa$ 
representation with value of the central charge given by the sum of the 
values (0 or 1, see discussion in Sec. \ref{sect3}) of the central charges 
associated to each line of the two-dimensional lattice. 

\medskip

In the previous sections we have discussed the case of $\vert q \vert
= 1$. The case of $q$ real can also be discussed and we refer to
\cite{LS93} for the definition of anyons for generic $q$.

\medskip

One can naturally ask if the realization in terms of $a$-anyons and
$A$-anyons here presented can be used to realize the deformation of
other finite or affine superalgebras. It seems that this procedure
can be extended to the other basic finite superalgebras, i.e. the
series $B(0, N)$, $B(M,N)$, $C(N+1)$, $D(M,N)$, while it is not clear
its extension to the exceptional or strange finite superalgebras or
to the affine case.
Finally we want briefly to comment on the difference between our
approach and the approach of \cite{LMR96}, even if here we present the
realization of a $q$-superalgebra and in \cite{LMR96} a realization of
a $q$-algebra is presented.  The approach of \cite{LMR96} is made on the
continuum and the authors do not use, as already remarked, $q$-bosons, 
as in the present paper, but standard bosons before ``anyonization''.
However one has to stress that their approach is based on a Fock
space realization which guarantees the consistency of the commutation
relations. It is worth noticing that on the Fock space fundamental
representation of the deformed algebra is indistinguishable from the
fundamental representation of the undeformed algebra. It follows that
a sum with the correct product of the fundamental representation
gives a representation of the deformed algebra. On the contrary
as we fulfill the commutation relations in abstract way, we are not
allowed to consider only the fundamental representation of the deformed
algebra realized by bosons and in order to achieve the consistency of the
representation we are lead to use $q$-bosons.

\section*{Acknowledgements}

Work supported by the European Commission TMR programme ERBFMRX-CT96-0045.
A.S. thanks the Laboratoire de Physique Th\'eorique ENSLAPP for kind 
hospitality during the period in which this paper was finished.

\newpage


\end{document}